\documentclass[12pt]{article}
\usepackage{amsfonts}
\begin{document}
\def\R{\mathbb{R}}
\def\C{\mathbb{C}}
\def\Z{\mathbb{Z}}
\def\N{\mathbb{N}}
\def\Q{\mathbb{Q}}
\def\D{\mathbb{D}}
\def\Sp{\mathbb{S}}
\def\T{\mathbb{T}}
\def\hb{\hfil \break}
\def\ni{\noindent}
\def\i{\indent}
\def\a{\alpha}
\def\b{\beta}
\def\e{\epsilon}
\def\d{\delta}
\def\De{\Delta}
\def\g{\gamma}
\def\qq{\qquad}
\def\L{\Lambda}
\def\E{\cal E}
\def\G{\Gamma}
\def\F{\cal F}
\def\K{\cal K}
\def\A{\cal A}
\def\B{\cal B}
\def\M{\cal M}
\def\P{\cal P}
\def\Om{\Omega}
\def\om{\omega}
\def\s{\sigma}
\def\t{\theta}
\def\th{\theta}
\def\Th{\Theta}
\def\z{\zeta}
\def\p{\phi}
\def\P{\Phi}
\def\m{\mu}
\def\n{\nu}
\def\l{\lambda}
\def\Si{\Sigma}
\def\q{\quad}
\def\qq{\qquad}
\def\half{\frac{1}{2}}
\def\hb{\hfil \break}
\def\half{\frac{1}{2}}
\def\pa{\partial}
\def\r{\rho}
\begin{center}
{\bf GAUSSIAN RANDOM FIELDS ON THE SPHERE AND SPHERE CROSS LINE}
\end{center}
\begin{center}
{\bf N. H. BINGHAM and Tasmin L. SYMONS}
\end{center}\begin{center}
In memory of Larry Shepp
\end{center}

\ni {\bf Abstract} \\
\i  We review the Dudley integral for the Belyaev dichotomy applied to Gaussian processes on spheres, and discuss the approximate (or restricted) continuity of paths in the discontinuous case.  We discuss also the spatio-temporal case, of sphere cross line.  In the continuous case, we investigate the link between the smoothness of paths and the decay rate of the angular power spectrum, following Tauberian work of the first author, Malyarenko, and Lang and Schwab.   \\

\ni {\bf Key words}.  Belyaev dichotomy, Dudley integral, Gaussian processes, spheres, ultraspherical polynomials, Schoenberg's theorem, Tauberian theorems, spherical functions, multiplication theorem, Feldheim-Vilenkin integral.  \\

\ni {\bf MSC Subject classification}. 60B15, 60B99, 60G15, 60G60.     \\ 

\ni {\bf 1.  Introduction} \\
\i Motivated by the mathematics of Planet Earth, we consider several aspects of Gaussian processes (Gaussian random fields) on spheres.  We begin (\S 2) with path-continuity: Belyaev's dichotomy and the Dudley integral, and the link between continuity and restricted measurability of the process. In \S 3 we combine the Karhunen-Lo\`eve expansion for a general Gaussian process with the spherical harmonics needed for the sphere, and deal also with the spatio-temporal case (sphere cross line).  In \S 4 we give a short proof of Malyarenko's theorem [Mal1,2], by Tauberian methods deriving from [Bin5,6], linking the asymptotics of ultraspherical series to the spherical case of the Dudley integral, and treat the infinite-dimensional case (the Hilbert sphere).  We also complement this by using Tauberian methods from [Bin6] to study integrability rather than asymptotics, as in the work of Lang and Schwab [LangS].     \\
\\
\\  

\ni {\bf 2.  Belyaev's dichotomy and the Dudley integral} \\

\i Let $X = \{ X_t: t \in M \}$ be a real-valued zero-mean Gaussian process, {\it on} (defined on, indexed by) $M$.  Here, $M$ will be the $d$-sphere ${\Sp}^d \subset {\R}^{d+1}$ of radius 1; the motivating example is $d = 2$, with $M = {\Sp}^2$ as Planet Earth.  The law of $X$ is determined by either of the covariance or the incremental variance:
$$
c(s,t) := cov(X_s, X_t) = E[X_s X_t], \qquad i(s,t) := E[(X_t - X_s)^2] 
$$      \\  
(respectively positive and negative definite, or of positive and negative type); we can pass between them by
$$
i(s,t) = c(s,s) + c(t,t) - 2 c(s,t), \quad
c(s,t) = \half (i(s,o) + i(t,o) - i(s,t)), 
$$
with $o$ some base point (North Pole).  We restrict attention to {\it isotropic} processes (those with {\it stationary increments}), where these are functions only of the geodesic
 distance $d(s,t)$, or of $x := \cos d(s,t) \in [-1,1]$ $(s,t \in M)$:
$$
c(s,t) = C(x), \qquad i(s,t) = I(x).
$$
We assume also that the covariance is {\it continuous}.  We can then use reproducing-kernel Hilbert spaces and the Karhunen-Lo\`eve expansion, which we will need below ([MarR, p.203-207], [Adl, III.2, III.3]). \\
\i We need the Gegenbauer (ultraspherical) polynomials $C_n^{\l}(x)$ [Sze, \S 4.7], normalised to take the value 1 at 1; for these we use Bochner's notation $W_n^{\l}(x)$ (see e.g. [Bin2]): 
$$
W_n^{\l}(x) := C_n^{\l}(x)/C_n^{\l}(1) 
= C_n^{\l}(x).n! \Gamma(2 \l)/\Gamma(n + 2 \l).
$$
These are the orthogonal polynomials generated by the probability measure
$$
G_{\l}(dx) := \frac{\G(\l + 1)}{\sqrt{\pi} \G(\l + \half)}.
(1 - x^2)^{\l - \half} dx \qquad (x \in [-1,1]):
$$
$$
\int_{-1}^1 W_m^{\l}(x) W_n^{\l}(x) G_{\l}(dx) 
= {\delta}_{mn}/{\om}_n^{\l},
\qquad
{\om}_n^{\l} 
:= \frac{(n + \l)}{\l}.\frac{\G(n + 2 \l)}{\G(2 \l)}.
$$
\i Half-integer values of the Gegenbauer index $\l$ correspond to (integer) values of the Euclidean dimension $d$ as above by
$$
\l = \half (d - 1);
$$
we will also need the Hilbert-space case $\l = \infty$) [Bin4]. \\ 
\i From the Bochner-Schoenberg theorem of 1940-42 ([Sch]; see [BinS1] for further references), $C$ is then, to within a scale factor $c \in (0,\infty)$ (reflecting physical units), a {\it mixture} (with mixing law $a = (a_n)_{n=0}^{\infty}, a_n \geq 0, \sum a_n = 1$) of {\it ultraspherical polynomials} $W_n^{\l}(x)$ with $\l := \half (d-1)$:
$$
C(x) = c \sum_0^{\infty} a_n W_n^{\l}(x), \qquad I(x) 
= c \sum_0^{\infty} a_n (1 - W_n^{\l}(x)).    \eqno (BS)
$$
\i The behaviour of the process is governed by the decay of the sequence $a = (a_n)$ above: the faster the decay, the better the behaviour (the smoother the paths, etc.)  The term {\it angular power spectrum (APT)} is used ([MariP]: see below) for a variant on $(a_n)$; we shall use the term here too, for convenience. \\
\i We recall Belyaev's dichotomy for Gaussian processes ([Bel]; [MarR, Th. 5.3.10]): Gaussian paths are either very nice (continuous), or very nasty (pathological: unbounded above and below on any interval, or set of positive measure).  Much is known, by way of necessary conditions for continuity [MarR, \S 6.2], and sufficient conditions [MarR, \S 6.1]; see e.g. [MarS1,2], [Gar].  One uses the {\it Dudley metric} (actually a pseudo-metric)
$$
d_X(s,t) := \sqrt{E[(X_s - X_t)^2]} \qquad (s, t \in M).
$$
For $u > 0$, write $N(u)$ for the minimum number of $d_X$-balls of radius $u$ needed to cover the parameter-space $M$; then if $H(u) := \log N(u)$, $H := \{ H(u): u > 0\}$ is called the {\it metric entropy}.  The {\it Dudley integral} is
$$
\int_0^{\e} \sqrt{H(u)} du \qquad (\e > 0).    \eqno(Dud)
$$ 
One obtains a clean necessary and sufficient condition for continuity, finiteness of $(Dud)$ [Dud1,2], only for $X$ isotropic [MarR], which is why we restrict to this here.  If $\phi$ is a non-negative function increasing near 0 with
$$
d_X(s,t) \leq \phi(|s - t|) \qquad (s, t \in M),
$$
then the Dudley integral is finite if
$$
\int_M^{\infty} \phi(e^{-x^2}) dx < \infty: \quad 
\int_0^{\e} \frac{\phi(u)}{\sqrt{- \log u}}.\frac{du}{u} < \infty.
$$  
For isotropic processes on spheres, take
$$
\phi(u) := \sup \{ \sqrt{I(cos v)}: v \leq u \}:
$$
the condition for path continuity of $X$ becomes ([Gar]; [Dud2, \S 7])
$$
\int_0^1 \sqrt{\frac{\sup_{v \leq u}I(\cos v)}{- \log u}} \frac{du}{u} < \infty: 
$$
$$
\int_0^1 \sqrt{\frac{\sup_{v \leq u} (1 - \sum_0^{\infty} a_n W_n^{\l}(\cos v))}
{- \log u}} \frac{du}{u} < \infty. \eqno(DudSph)
$$ 
\i Despite the `pathological' behaviour of the sample paths of the process in the discontinuous case of Belyaev's dichotomy, there is a sense in which they are `nearly continuous': a `localisation of pathology'.  For, by the Karhunen-Lo\`eve expansion,
$$
X(t, \om) 
= \sum_0^{\infty} {\phi}_n(t) Z_n(\om),        \eqno(KL)
$$
with the $Z_n$ independent standard normal random variables and the ${\phi}_n$ continuous functions (a.s.: see [MarR, Remark 5.3.3]).  In particular, a.s. (we may exclude the exceptional $P$-null $\om$-set from our sample space and so omit this restriction), $X(t) = X(t, \om)$ is a {\it measurable} function of $t$.  So, by Lusin's continuity theorem (or Lusin's restriction theorem, of 1912: [Dud3, Th. 7.5.2], [Rud, \S 2.24]), $X(t)$ becomes {\it continuous} in $t$ when restricted to a time-set of $t$ avoiding a set of arbitrarily small measure. \\

\ni {\bf Remarks}. \\
1. {\it The oscillation function}. \\
\i The local behaviour of the paths of $X(t, \om)$ is governed by the {\it oscillation function}, a deterministic function, $\a(t)$ say [MarR, p.209-211].  The continuous case of Belyaev's dichotomy has $\a \equiv 0$, the discontinuous case has $\a \equiv \infty$.  Despite the `pathological' appearance of this case, it is as well to note that a measurable function can have oscillation function $\equiv \infty$, as here. \\
2. {\it Approximate limits and limsups}. \\
\i The Lusin argument above from $(KL)$ says that the paths, when slightly restricted, become continuous in some sense.  The concepts of {\it approximate limit}, $ap-lim$ (and so of approximate limsup, $ap-limsup$, and approximate derivative) and {\it approximate continuity} go back to Denjoy in 1916 and Khintchine in 1927 (see e.g. Saks [Sak, IX.10]).  An approximate limit at a point $t$ becomes an actual limit when the approach to $t$ is made through a Borel set having $t$ as a {\it density point} (in the sense of the Lebesgue density theorem, see e.g. [Rud, \S 7.2], [Dud3, p.422], [BinO2, Th. L], and of the density topology, below).  For a number of probabilistic results on $ap-limsup$, see Geman and Horowitz [GemH, \S 13] (cf. [Adl, IV.6]).  For the relevant real-variable theory, see e.g. [GemH, \S 14], and the earlier paper by Smallwood [Sma]. \\
3. {\it The density topology}. \\
\i The density topology takes as its open sets the measurable sets all of whose points are density points.  That this gives a topology, and the intimate link with Denjoy's approximate continuity (and so with $ap-lim$), are due to Haupt and Pauc [HauP].  It has been much studied, by C. Goffman and others; for references see e.g. [BinO1,3]. \\

\ni {\bf 3. The Karhunen-Lo\`eve expansion; the spatio-temporal case} \\

\ni {\it Spherical harmonics} \\ 
\i For Gaussian processes on (parametrised by) spheres, as here, one appropriate system to use for spectral expansions as in $(KL)$ is the {\it spherical harmonics} ([AndAR, Ch. 9]; [SteW, IV.2]).  These are the restrictions to the $d$-sphere ${\Sp}^d \subset {\R}^{d+1}$ of homogeneous {\it harmonic} polynomials -- solutions to Laplace's equation in ${\R}^{d+1}$.  For each degree $\ell = 0,1,2,\cdots$, there are
$$
c(\ell,d) := \frac{2\ell + d - 1}{d - 1} {\ell + d - 2 \choose \ell}
$$
linearly independent spherical harmonics of degree $\ell$.  Furthermore, the usual time parameter of a stochastic process is replaced by a space parameter, $x \in {\Sp}^d$.  Using now $T = T(x)$ for the process, now a {\it random field}, it has a KL expansion of the form
$$
T(x) 
= \sum_{\ell = 0}^{\infty} \sum_{m = 0}^{c(\ell,d)}
a_{\ell m} Y_{\ell m}(x)
= \sum_{\ell m} a_{\ell m} Y_{\ell m}(x),
= \sum_{\ell m} \langle T, Y_{\ell, m} \rangle Y_{\ell m}(x),
$$
say.  Here the randomness is in the {\it random Fourier coefficients}
$$
a_{\ell m} = \langle T, Y_{\ell, m} \rangle.
$$  
By $(KL)$, they are independent if the process $T$ is Gaussian, and conversely ([LusP], above).  Baldi and Marinucci [BalM] derive the converse from the Skitovich-Darmois theorem, one of many {\it characterization theorems} for normality (Gaussianity); see e.g. 
[KagLR].  \\

\ni {\it Spherical harmonics and Gegenbauer polynomials}\\
\i The link between the spherical harmonics $Y_{\ell m}$ above and the (normalized) Gegenbauer polynomials $W_n$ is the {\it addition theorem for spherical harmonics} ([AndAR, Th. 9.6.3]; [Yad, I.5.1]):
$$ 
\sum_{m=1}^{c(\ell, d)} Y_{\ell m}(x) Y_{\ell m}(y) 
= W_{\ell}^{\l}(\cos \langle x,y \rangle).
c(\ell, d)/{\om}_d,                           \eqno (Add)
$$
with $\langle x, y \rangle$ the inner product of $x$ and $y$ in ${\R}^{d+1}$ (so $\cos \langle x, y \rangle$ is the geodesic distance between $x$ and $y$ in ${\Sp}^d$), and ${\om}_d$ the surface area of ${\Sp}^d$. \\
\i The isotropy of $T$ is reflected in that the laws of the $a_{\ell m}$ depend only on the degree $\ell$ of the spherical harmonic and not on $m$:
$$
a_{\ell m} \sim N(0, v_l),
$$
say (`$v$ for variance').  The covariance is calculated by  
\begin{eqnarray*}
C(x,y) := cov(T(x), T(y)) 
&=& E[(\sum_{\ell m} a_{\ell m} Y_{\ell m}(x))
(\sum_{{\ell}' m'} a_{{\ell}' m'} Y_{{\ell}' m'}(y))] \\
&=& \sum_{\ell m} \sum_{{\ell}' m'} E[a_{\ell m} a_{{\ell}' m'}] 
Y_{\ell m}(x) Y_{{\ell}' m'}(y)] \\
&=& \sum_{\ell m} \sum_{{\ell}' m'} {\delta}_{\ell {\ell}'}{\delta}_{{\ell}' m'} v_{\ell} Y_{\ell m}(x) Y_{{\ell}' m'}(y)] \\
&=& \sum_{\ell} v_{\ell} \sum_m Y_{\ell m}(x) Y_{\ell m}(y): 
\end{eqnarray*}
$$
C(x,y)
= \frac{1}{\om_d}\sum_{\ell} v_{\ell} c(\ell, d)
W_{\ell}^{\l}(\cos \langle x,y \rangle),
$$ 
which (to within notation) is Schoenberg's theorem.  We write
$$
a_{\ell} := v_{\ell} c(\ell,d)/{\omega}_{\ell},
$$
and call $a = (a_{\ell})$ the {\it angular power spectrum} (APS).  The covariance and the KL expansion then become
$$
C(x,y)
= \sum_{\ell} a_{\ell} W_{\ell}^{\l}(\cos \langle x,y \rangle),
                                              \eqno (APS)
$$
$$
T(x) = \sum_{\ell m} a_{\ell m} Y_{\ell m}(x), \q a_{\ell m} \sim N(0,v_{\ell}),                       \eqno (T(x))
$$
with the $a_{\ell m}$ independent. \\
\ni {\it Note}.  Some authors (e.g. [MariP]) use the term angular power spectrum for $v = (v_{\ell})$ instead.  For the case $d = 2$, conversion between the two becomes especially simple: $c(\ell, 2) = 2 \ell + 1$, so
$$
a_{\ell} := v_{\ell}. (2 \ell + 1)/4 \pi.
$$

\ni {\it Spatio-temporal version} \\
\i The covariance calculation above can be extended to the spatio-temporal setting.  We already know the form of the general covariance, from the Berg-Porcu theorem [BerP] (see $(BP)$ below).  The spatio-temporal covariances here have a very similar structure to the Bochner-Schoenberg form $(BS)$ in the spatial case,
$$
c \sum_0^{\infty} a_n W_n^{\lambda}(\cos \langle x, y \rangle), 
\quad c > 0, \ a_n \geq 0, \ \sum a_n = 1       \eqno (BS)
$$
(to within a constant multiple, a mixture of ultraspherical polynomials), but now as
$$
c \sum_0^{\infty} a_n W_n^{\lambda}(\cos \langle x, y \rangle) {\phi}_n(t - s), 
\quad c > 0, \ a_n \geq 0, \ \sum a_n = 1,      \eqno(BP)
$$
with the ${\phi}_n$ characteristic functions of probability laws on the line.  The role of the {\it Bochner-Schoenberg theorem} in $(BS)$ is played in $(BP)$ by the {\it Bochner-Godement theorem} ([AskB], [BinS1]; the result is also derived in [Ma2]).  Many results proved for the class $\mathcal{P}(\Sp^d)$ thus carry over to $\mathcal{P}(\Sp^d \times \R)$.  So we need a choice of expansion $(KL)$ that will lead to this.  We generalise $(T(x))$ above in the result below.  Note that for a centred Gaussian process, the covariance is canonical (it determines the process), but the Karhunen-Lo\`eve expansion is not (it depends on the complete orthonormal system).  What one needs is a KL expansion that generates the required covariance. \\

\ni {\bf Theorem (Karhunen-Lo\`eve expansion for sphere cross line}).  The Karhunen-Lo\`eve expansion
$$
T(x,t) = \sum_{\ell m} a_{\ell m} c_{\ell}{\phi}_{\ell}(t) Y_{\ell m}(x), \q 
a_{\ell m} \sim N(0,v_{\ell}),             \eqno (T(x,t))
$$
with $Y_{\ell m}(x)$ the spherical harmonics, ${\phi}_{\ell}(t)$ a characteristic function and $c_{\ell}$ constants with $\sum_{\ell} v_{\ell} c(\ell, d) c_{\ell}^2 < \infty$, has covariance the Berg-Porcu form $(BP)$ above.  It thus generates the most general Gaussian random field on ${\Sp}^d \times \R$ isotropic on ${\Sp}^d$ and stationary on $\R$. \\ 

\ni {\it Proof}.  The covariance calculation now becomes
\begin{eqnarray*}
C((x,s), (y, t)) 
&=& cov(T(x, s), T(y,t))  \\ 
&=& E[(\sum_{\ell m} a_{\ell m}c_{\ell} {\phi}_{\ell}(s)  Y_{\ell m}(x))
(\sum_{{\ell}' m'} a_{{\ell}' m'} c_{{\ell}'} {\phi}_{{\ell}'}(t) Y_{{\ell}' m'}(y))] \\
&=& \sum_{\ell m} \sum_{{\ell}' m'} 
E[a_{\ell m} a_{{\ell}' m'}]
c_{\ell} c_{{\ell}'} {\phi}_{\ell}(s) {\phi}_{\ell}(t) 
Y_{\ell m}(x) Y_{{\ell}' m'}(y) \\
&=& \sum_{\ell m} \sum_{{\ell}' m'} 
{\delta}_{\ell {\ell}'}{\delta}_{{\ell}' m'} v_{\ell}
c_{\ell} c_{{\ell}'} {\phi}_{\ell}(s) {\phi}_{\ell}(t)  Y_{\ell m}(x) Y_{{\ell}' m'}(y) \\
&=& \sum_{\ell} v_{\ell} c_{\ell}^2 {\phi}_{\ell}(s) {\phi}_{\ell}(t) \sum_m Y_{\ell m}(x) Y_{\ell m}(y):
\end{eqnarray*}
$$
C((x,s), (y, t)) 
= 
\frac{1}{\om_d}\sum_{\ell}v_{\ell}c(\ell,d)c_{\ell}^2.
{\phi}_{\ell}(s) {\phi}_{\ell}(t) W_{\ell}(\langle x,y \rangle),
$$
by $(Add)$. \\
\i We are assuming isotropy in time (stationarity) as well as in space.  So the covariance depends only on the difference of the two times $s$ and $t$.  So we may replace $s$ by  $0$ (so ${\phi}_n(0) = 1$) and $t$ by $t-s$.  This now reduces to $(BP)$ above, the Berg-Porcu form [BerP].  As we know that this form for the covariance is the most general one, and a (centred) Gaussian process (random field here) is specified by its covariance, this shows that no generality has been lost by the above choice for the Karhunen-Lo\`eve expansion.   \hfil $\square$ \break

\ni {\it Note}. 1.  The earliest work of this kind that we are aware of is the pioneering paper of Jones [Jon] in 1963.  This was recently extended from ${\Sp}^2$ to ${\Sp}^d$ by Clarke De la Cerda, Alegria and Porcu [ClaAP].  Both approaches deal with the space and time aspects asymmetrically, and neither has the randomness in the simplest form possible: a {\it sequence of independent standard normals}, as in $(KL)$ above (we recommend [MarR \S 5.3] for a treatment of Karhunen-Lo\`eve expansions).  Our approach above is both symmetrical and simpler; this enables us to avoid the difficulties touched on at the end of [ClaAP, \S 3.1]. \\
\ni 2. In [ClaAP], distance on sphere cross line is taken under the cartesian rule $s^2 = s_1^2 + s_2^2$ for a product of metric spaces, whereas above we use the differential cartesian rule $ds^2 = ds_1^2 + ds_2^2$ for a product manifold.  The difference is not important here, but can be crucial.  For an instance, see [BinMS]: {\it  Brownian motion exists} on sphere cross line in the first case (and under Hamming distance $s = s_1 + s_2$), but not the second. \\
            
\ni {\bf 4.  Malyarenko's theorem; the Hilbert sphere}\\

\i What $(DudSph)$ above says is that the paths are continuous if and only if the coefficients $a_n$ in the mixture law $a$ in $(BS)$ -- the {\it angular power spectrum} -- {\it decay fast enough}.  Slow decay means wild behaviour of the paths, but if the decay is fast enough, the paths become very smooth.  As we shall see, if $a_n = O(1/n^{1+\a})$ for $\a > 0$, the paths are continuous (and become smoother with increasing $\a$).   \\
\i While the condition $(DudSph)$ resolves the matter completely in principle, in practice implementing it is formidable, for the obvious three reasons: passage between the mixing law $(a_n)$ and the ultraspherical series $\sum a_n W_n^{\l}$, the supremum, and the integration.  The nub here is the first: the link between the {\it decay} of $a_n$ for large $n$, and the {\it growth} of $1 - \sum_0^{\infty} a_n W_n^{\l}(\cos v)$ for small $v > 0$. \\
\i While there is no definitive answer to this question (any more than there is in the classical case of Fourier series [Zyg]), there is an answer in the principal case of practical interest, that when the angular power spectrum $(a_n)$ is {\it regularly varying} (see e.g. [BinGT]).  Here the results (whose proofs we sketch below) are due to Malyarenko [Mal1,2], based on early work of the first author [Bin5] (itself based on earlier work of Askey and Wainger [AskW]): \\

\ni {\bf Theorem (Malyarenko)}.  For $\ell$ slowly varying,
$$
A_n := \sum_n^{\infty} a_k \sim \ell(n)/n^{\g} \quad (n \to \infty) \quad (\g \in (0,2))
$$
iff
$$
I(v) = 1 - \sum_0^{\infty} a_n W_n^{\l}(\cos v) \sim 
\G(\l + \half).
\frac{\G(1 - \half \g) }{2^{\g} \G(\l + \half - \half \g)}. 
v^{\g} \ell(1/v) \quad (v \downarrow 0).
$$
\ni {\it Proof}. The implication from $A_n$ to $I(v)$ is Abelian; the converse is Tauberian.  One has
$$
I(v) = \sum_0^{\infty} a_n (1 - W_n^{\l}(v)
= \sum_0^{\infty} (A_n - A_{n+1})(1 - W_n^{\l}(v)),
$$
and writes this by partial summation as
$$
I(v) = \sum A_{n+1}(W_n^{\l}(v) - W_{n+1}^{\l}(v)). 
$$         
\i The difference of ultraspherical polynomials here may be expressed as {\it one} Jacobi polynomial (Erd\'elyi et al. [ErdMOT, Vol. II, 10.8(32)].  Recall that the Jacobi polynomials are a two-index family $P_n^{(\a, \b)}$   ($\a, \b \geq - \half$; we take $\a \geq \b$).  When $\a = \b$, the Jacobi polynomials reduce to the ultraspherical polynomials, with (as above)
$$
\a = \b = \l - \half = \half(d - 2).
$$
We use the normalisation [Mal2, 4.3.1]
$$
R^{(\a, \b)}_n(x) := P^{(\a, \b)}_n(x)/P^{(\a, \b)}_n(1).
$$  
Then ([Mal2, p.127], [ErdMOT II 10.8(32)])
$$
R_n^{\a, \b}(\cos \th) - R_{n+1}^{\a, \b}(\cos \th) 
= \frac{(2n + \a + \b + 2)}{(\a + 1)} {\sin}^2 \half \th \ R_n^{\a + 1, \b}(\cos \th).
$$
So
$$
I(\cos \th) = \frac{2{\sin}^2 \half \th }{(\a + 1)}  
\sum (n + \a + 1) A_n 
R_n^{\a + 1, \b}(\cos \th).                  \eqno(\ast)
$$  
The ${\sin}^2 \half \th$ (equivalently, ${\th}^2/4$) term on the right of $(\ast)$ accounts for the upper limit 2 on $\g$ in the result; that the incremental variance is non-negative accounts for the lower limit of 0.  The results of [Bin5] now apply to the sequence $(n + \a + 1)A_n = (n + \a + 1) \sum_n^{\infty} a_k$ with the $\sigma$ there as $1 - \g$.   The Tauberian conditions needed follow from $a_n \geq 0$ (so $A_n$ is non-negative and non-decreasing). \hfil $\square$ \break

\i Malyarenko's theorem is very similar to that of [Bin3] (proved more simply in [Bin7]) on Hankel transforms, the link being provided by Szeg\H{o}'s Hilb-type asymptotic formula for the Jacobi polynomials [Sze, Th. 8.21.12]. \\
\i The Belyaev integral is of course convergent in all these cases, and so Malyarenko's theorem provides us with an ample range of examples of the continuous case in the Belyaev dichotomy (the pathological case being of course less common in practice).  For, the supremum operation in $(DudSph)$ (a reflection of the great mathematical difficulties in bridging the gap between the necessary and the sufficient conditions for finiteness of the Dudley integral) is harmless here: any regularly varying function of {\it non-zero} index is asymptotically monotone [BinGT, \S 1.5.2]).  \\
\i One can extend to $\g = 0$ here, when the tail $A_n$ of the mixing law is slowly varying, but convergence of the Dudley integral now hinges on the behaviour of $\ell$ at infinity.  This is shown by familiar examples such as $\sum 1/(n (\log n)^k)$, $\sum 1/(n \log n (\log \log n)^k)$, each convergent if $k > 1$, divergent if $k \leq 1$.  One can also extend to the case $\g = 2$ [Bin5].  \\
\i $O$-versions of these results are straightforward (cf. Korevaar [Kor, IV.10]).  \\

\ni {\it The Hilbert sphere}. \\
\i The Hilbert sphere ${\Sp}^{\infty}$ is not locally compact, and because of this one may expect very different behaviour for it from that on Euclidean spheres.  Gaussian processes on ${\Sp}^{\infty}$ are {\it discontinuous} (L\'evy [Lev]; Berman [Berm1,2]; Dudley [Dud2, \S 5]).  More is true: such processes are {\it locally deterministic} (see [Lev, p.355], [Berm1, p.950] for the definition): the behaviour of the process locally determines it everywhere.  This sounds reminiscent of the great smoothness shown by holomorphic functions in complex analysis, but is in fact diametrically opposite: the process is extremely wild, and `gets everywhere it will go immediately'. \\
\i The ultraspherical polynomials may be defined for  $\l = \infty$ by $W_n^{\infty}(x) = x^n$ (see e.g. [Bin1]).  But, as
$$
\Gamma ( \l + \half + \half \gamma ) / \Gamma ( \l + \half) 
\sim {\l}^{\half \gamma } \to \infty \quad (\l \to \infty),
$$
this case does not follow formally from Malyarenko's theorem by letting $\l \to \infty$.  Instead, we have here: \\

\ni {\bf Proposition}.  In the notation of Malyarenko's theorem,
$$
A_n := \sum_n^{\infty} a_k 
\sim \ell(n)/n^{\half \g} \quad 
(n \to \infty) \quad (\g \in (0,2))
$$
iff
$$
1 - \sum a_n (\cos v)^n \sim \frac{\G (1 - \half \gamma )}{2^{\half \gamma}} v^\l \ell(1/v^2) \quad (v \downarrow 0).
$$

\ni {\it Proof}.  The functions in ${\cal P}_{\infty}$ are the probability generating functions (in $t$, say), or (putting $t =
e^{-s}$) the Laplace-Stieltjes transforms.  We can read off the relevant tail-behavour here from e.g. [BinGT, Cor. 8.1.7]. Writing $\cos v = e^{-s}$ here, we have $s \sim \half v^2$ as $s, v \downarrow 0$, which gives the result.  \hfil $\square$ \break

\i \i  One must expect the tails in the Hilbert case here (with the sphere non-compact) to be heavier than in the Euclidean case of Malyarenko's theorem (with the sphere compact): now, the paths are wild rather than continuous, and there there are `more ways of going off to infinity'.  Thus the relevant probability laws $(a_n)$ here have regularly varying tails with index in $(0,1)$, rather than in $(0,2)$ as Malyarenko's theorem -- that is, they correspond to {\it infinite mean} rather than {\it infinite variance}. \\
\i The constants introduced (in going between the `Abelian' and `Tauberian' sides) in results of this type are the values, for $s = \g$, of the {\it Mellin transform}
$$
\hat k(s) := \int_0^{\infty} u^s k(u) du/u \quad (s \in \C)
$$
of the kernel $k$ in the relevant Mellin-Stieltjes convolution (see e.g. [BinGT, Ch. 4, 5]).  In the Hankel case of [Bin3, 7] the relevant transform is exactly of convolution type; here and in [Bin5] the `ultraspherical transform' is only approximately so (cf. [BinGT, \S 4.2, 4.3, 4.10]).  It is interesting to compare the Mellin transforms in these three cases. \\  

\ni {\bf Remarks}. \\
\ni 1. {\it Besov paths}. \\
\i Kerkyacharian et al. [KerOPP, \S 7.22] show that if the angular power spectrum satisfies $a_n = O(1/n^{1 + \g})$ for $\g > 0$ (so $A_n := \sum_n^{\infty} a_k = O(1/n^{\g})$), then the sample paths of the process $X$ are a.s. in the Besov space $B^{\a}_{\infty, 1}$ for all $\a < \g$ (see Gin\'e and Nickl [GinN] for the theory of Besov spaces in such contexts, Fukushima et al. [FukOT] for the necessary theory of Dirichlet structure on the index set, ${\Sp}^d$ here).  Thus {\it the faster the decay of the angular power spectrum, the smoother the paths of the process}. \\
\ni 2. {\it Fractional calculus on spheres}. \\
\i Following Askey and Wainger [AskW, Part I Section III], a theory of fractional integration and differentiation on spheres was given by Bavinck [Bav].  This is based on the expansion into spherical harmonics $S_{l,m}$ above; these are eigenfunctions of the spherical Laplacian $\De$ (Laplace-Beltrami operator on the sphere), with eigenvalues $-l(l+\l)$ (or $-l(l + \a + \b + 1)$ in the Jacobi case):
$$
\De S_{l,m} = -l(l+\l)S_{l,m}, \quad
(1 - \De)S_{l,m} = (1 + l(l + \l))S_{l,m}.
$$  
In terms of the {\it fractional Laplacian} (see e.g. [Hor], [Stei, V.1.2]), applying $(1 - \De)^{\s/2}$ introduces multipliers $(1 + l(l+\l))^{\s/2}$ into the expansion.  For $\s > 0$, this corresponds to (fractional) differentiation of order $\s$ ($\De$ being a second-order differential operator), or (fractional) integration if $\s$ is negative (recall: the faster the angular power spectrum coefficients decay, the smoother the paths of the process, and the slower, the rougher).        \\
\i This has the semigroup property
$$
I_{\a + \b} = I_{\a} \circ I_{\b}.
$$
This desirable property is not shared by previous definitions of spherical fractional integration (see [AskW] for references), nor by analogues in the literature on `dimension walks'; see e.g. [BinS3] for references. \\
\i Note that $\l$ here may be a {\it continuous} parameter, and is not restricted to the half-integer values implied by $\l = \half (d - 1)$ with $d$ the dimension of the sphere (as a Riemannian manifold).  See [Bin1], [BinS3] for projections between two different dimensions (parameters). \\

\ni {\bf 5.  Integrability and path-continuity} \\

\i The question of path-continuity of the process is addressed in the work of Lang and Schwab [LangS] (cf. [AndL]) and Lan, Marinucci and Xiao [LanMX].  The picture is much as above: the faster the decay of the angular power spectrum, the better: the more regular the paths of the process (and, as in Malyarenko's theorem of \S 3,  the faster the decay of the incremental variance at the origin). \\
\i In [LangS, \S 4, Assumption 4.1], Lang and Schwab assume a decay condition on the angular power spectrum measured by a {\it summability condition} (rather than by rate of decay as in Malyarenko's theorem): in our notation, they assume
$$
\sum a_n n^{\g} < \infty \qquad (\g > 0).           \eqno(Int)                         
$$
In view of $(\ast)$ above, we re-write this by partial summation as
$$
\sum A_n. n^{\g - 1} < \infty: \qquad 
\sum (n+\a + 1)A_n.n^{\g - 2} < \infty.
$$
\i As above, and in [LangS \S 4], the case $\g \in (0,2)$ is specially important, so we begin with that. Then the summability condition $(Int)$ may [Bin6, Th. 1] be translated into a corresponding {\it integrability condition} on the incremental variance at the origin: $(Int)$ implies
$$
\int_{0+}^{\pi/2} I(\cos \t).{\t}^{- \g} d\t/\t 
< \infty.                                             \eqno(Int')
$$ 
As $\int_{0+} d\t/\t$ diverges, this gives in particular that
$$
I(\cos \t) = o({\t}^{\g}) \quad (\t \downarrow 0).
$$
This strengthens the result of [LangS, Lemma 2] from $O(.)$ to $o(.)$ (though in view of the `$\epsilon$-gap' in [LangS, Th. 4.7], where it is used, this does not matter). \\
\i This leads quickly to the path-regularity result ([LangS]; cf. [LanMX]): \\

\ni {\bf Theorem (Lang and Schwab, [LangS Th. 4.7]}.  Under the summability condition $(Int)$ on the angular power spectrum, for any $\d < \g/2$ the process has a $C^{\d}$-valued modification:for $k$ the integer part of $\g/2$, the modification is $k$ times continuously differentiable, with $k$th derivative H\"older continuous with exponent $\d - k$. \\

\i The proof involves the following: \\
(i) For $n \in N$, $x, y \in {\Sp}^d$
$$
E[|X(x) - X(y)|^{2n}] \leq C_{\g,n} d(x,y)^{\g n},
$$
with $d(.,.)$ geodesic distance as before [LangS, Lemma 4.3]. \\
(ii) The Kolmogorov-Chentsov theorem on manifolds [AndL] gives the result for $\g \in (0,2]$. \\
(iii) For $\g > 2$, $k$-fold fractional differentiation (see \S 3 Remark 2) reduces to the range above.  \\
We refer for detail to [LangS], [AndL].  We return to such matters elsewhere. \\

\ni {\bf Remarks}. \\
\ni 1. {\it Strong local non-determinism}. \\
\i Using the concept of strong local determinism, Lan, Marinucci and Xiao [LanMX] improve the Lang-Schwab result above, obtaining an exact modulus of continuity (and so avoiding an $\epsilon$-gap), for the case $d = 2$ and with the angular power spectrum coefficients bounded above and below by constant multiples of powers.  This condition holds, for example, for spherical fractional Brownian motion (Lan and Xiao [LanX]). \\
2.  {\it Vector data}. \\
\i Often data on spheres are vectors, as several different quantities are measured (temperature, wind speed, humidity etc.); the relevant covariances are then matrices.  See e.g. [Ma1], where a number of applications are given. \\
3.  {\it Statistics}. \\
\i One extremely important application for the theory of Gaussian random fields on spheres is of course the study of cosmic microwave background (CMB) radiation; see [MariP] for a monograph treatment.  For statistical estimation in this and related areas, see e.g. Durastanti, Lan and Marinucci [DurLM], Leonenko, Taqqu and Terdik [LeoTT]. \\
4.  {\it Stochastic partial differential equations (SPDEs)}. \\
\i For the stochastic heat equation on the sphere, see Lang and Schwab [LangS, \S 7]. \\ 
\ni 5. {\it Regular variation and function spaces}. \\
\i To avoid the `$\epsilon$-gap' in the Lang-Schwab theorem above, one needs a finer scale of spaces than is provided by the powers (in particular, one that is sensitive to logarithmic factors, etc.)  One such is provided by the {\it Orlicz spaces}; see e.g. Krasnoselkii and Rutickii [KraR].  These lead to the {\it Besov-Orlicz spaces}; see e.g. [CieKR]. \\
\ni 6. {\it Integrability theorems for Fourier series}. \\
\i Much is known about integrability conditions for Fourier series.  For detail, see the two monographs on the subject, by Boas [Boa] and Yong [Yon] (as well as [Bin6] and the references cited there). \\

\ni {\bf Acknowledgements}.  The second author acknowledges financial support from the EPSRC Centre for Doctoral Training in the Mathematics of Planet Earth [EP/L016613/1].  Both authors thank Adam Ostaszewski for discussions, and the editors for the invitation to contribute to the Larry Shepp Memorial Issue.\\              

\ni {\bf Postscript}.  It is a pleasure for the first author to record here his happy memories of all his dealings with Larry Shepp, and of the excellent conference in his memory at Rice University, 25-29 June 2018, so ably organised by Philip Ernst.  The title of his talk there was `Four themes from the work of Larry Shepp', the first of which was Gaussian processes, as here. \\
\i It is also a pleasure for both authors to record here the great debt that they and their colleagues owe to the work of Larry Shepp. \\

\begin{center}
{\bf References}
\end{center}
\ni [Adl] R. J. Adler, {\sl An introduction to continuity, extrema and related topics for general Gaussian processes}.  Inst. Math. Statist. Lecture Notes- Monographs {\bf 12}, IMS, 1990. \\
\ni [AndL] R. Andreev and A. Lang, Kolmogorov-Chentsov theorem and differentiability of random fields on manifolds.  {\sl Potential Analysis} {\bf 41} (2014), 761-769. \\
\ni [AndAR] G. E. Andrews, R. Askey and R. Roy, {\sl Special functions}.  Cambridge University Press, 1999. \\
\ni [AskB] Askey, R. A. and Bingham, N. H.: Gaussian processes on compact symmetric spaces.  {\sl Z. Wahrschein.  verw. Geb.} {\bf 37} (1976), 127 -- 143. \\       
 \ni [AskW] R. Askey and S. Wainger, On the behaviour of special classes of ultraspherical polynomials, I, II.  {\sl J. Anal. Math.} {\bf 15} (1965), 193-220, 221-244. \\
 \ni [BalM] P. Baldi and D. Marinucci, Some characterizations of the spherical harmonics coefficients for isotropic random fields.  {\sl Statist. Probab. Letters} {\bf 77} (2007), 490-496. \\
\ni [Bav] H. Bavinck, A special class of Jacobi series and some applications.  {\sl J. Math. Anal. Appl.} {\bf 37} (1972), 767-797. \\
 \ni [Bel] Belyaev (Belayev), Yu. K., Continuity and H\"older's conditions for sample functions of stationary Gaussian processes.  {\sl Proc. Fourth Berkeley Symp. Math. Stat. Probab. Volume II: Contributions to Probability Theory} (ed. J. Neyman), 23-33, U. California Press, 1961. \\
\ni [BerP] C. Berg and E. Porcu, From Schoenberg coefficients to Schoenberg functions.  {\sl Constr. Approx.} {\bf 45} (2017), 217-241. \\ 
\ni [Berm1] S. M. Berman, A Gaussian paradox: determinism and discontinuity of sample paths.  {\sl Ann. Prob.} {\bf 2} (1974), 950-953. \\
\ni [Berm2] S. M. Berman, Isotropic Gaussian processes on the Hilbert sphere.  {\sl Ann. Prob.} {\bf 6} (1980), 1093-1106. \\ 
\ni [Bin1] N. H. Bingham, Integral representations for ultraspherical polynomials.  {\sl J. London Math. Soc.} {\bf 6} (1972), 1-11. \\
\ni [Bin2] N. H. Bingham, Random walk on spheres.  {\sl Z. Wahrscheinlichkeitstheorie  verw. Geb.} {\bf 22} (1972), 169-192.\\
\ni [Bin3] N. H. Bingham, Tauberian theorems for integral transforms of Hankel type.  {\sl J. London Math. Soc.} (2) {\bf 5} (1972), 493-503. \\
\ni [Bin4] N. H. Bingham, Positive definite functions on spheres.  {\sl Proc. Cambridge Phil. Soc.} {\bf 73} (1973), 145-156.  \\
\ni [Bin5] N. H. Bingham, Tauberian theorems for Jacobi series.  {\sl Proc. London Math. Soc.} {\bf 36} (1978), 285-309. \\
\ni [Bin6] N. H. Bingham, Integrability theorems for Jacobi series.  {\sl Publ. Inst. Math. Beograd} {\bf 26} (40) (1979), 45-56. \\
\ni [Bin7] N. H. Bingham, On a theorem of Klosowska about generalised convolution. {\sl Colloq. Math.} {\bf 48} (1984), 117-125. \\
\ni [BinGT] Bingham, N. H., Goldie, C. M. and Teugels, J. L., {\sl Regular variation}.  Cambridge University Press, 1987 (2nd ed. 1989). \\ 
\ni [BinMS] N. H. Bingham, A. Mijatovi\'c and Tasmin L. Symons, Brownian manifolds, negative type and geotemporal convariances.  {\sl Comm. Stoch. Analysis} (Herbert Heyer Festschrift) {\bf 10} no. 4 (2016), 421-432. \\
\ni [BinO1] N. H. Bingham and A. J. Ostaszewski, Beyond Lebesgue and Baire II: bitopology and measure-category duality.  {\sl Colloq. Math.} {\bf 121} (2010), 225-238.\\
\ni [BinO2] N. H. Bingham and A. J. Ostaszewski, Kingman, category and combinatorics. {\sl Probability and Mathematical Genetics} (Sir John Kingman Festschrift, ed.  N. H. Bingham and C. M. Goldie), 135-168, London Math.  Soc.  Lecture Notes in Mathematics {\bf 378}, Cambridge University Press, 2010.\\
\ni [BinO3] N. H. Bingham and A. J. Ostaszewski, Beyond Lebesgue and Baire IV: density topologies and a converse Steinhaus-Weil theorem.  {\sl Topol. Appl.} {\bf 239} (2018), 274-292. \\
\ni [BinS1] N. H. Bingham and Tasmin L. Symons, Dimension walks on ${\Sp}^d \times \R$.  {\sl Statistics and Probability Letters} {\bf 147} (2019), 12-17; arXiv1809.03955. \\
\ni [BinS2] N. H. Bingham and Tasmin L. Symons): Probability and statistics of Planet Earth. I: The Bochner-Godement theorem and geotemporal covariances. arXiv:1706.02972;  arXiv:1707.05205. \\
\ni [BinS3] N. H. Bingham and Tasmin L. Symons.  Integral representations for ultraspherical polynomials II.  Preprint. \\ 
\ni [Boa] R. P. Boas, {\sl Integrability theorems for trigonometric transforms}.  Erg. Math. {\bf 18}, Springer, 1967. \\
\ni [CieKR] Z. Ciesielski, G. Kerkyacharian and B. Roynette, Quelques espaces fonctionnels associ\'es \`a des processus Gaussiens. {\sl Studia Math.} {\bf 107} (1993), 171-214. \\
\ni [ClaAP] J. Clarke De la Cerda, A. Alegria and E. Porcu, Regularity properties and simulations of Gaussian random fields on the sphere cross time.  {\sl Electronic J. Stat.} {\bf 12} (2018), 399-426. \\
\ni [Dud1] Dudley, R. M., The sizes of compact subsets in Hilbert space and continuity of Gaussian processes.  {\sl J. Functional Analysis} {\bf 1} (1967), 290-330.   \\
\ni [Dud2] Dudley, R. M., Sample functions of the Gaussian process.  {\sl Ann. Probab.} {\bf 1} (1973), 66-103.   \\
\ni [Dud3] R. M. Dudley, {\sl Real analysis and probability}.  Wadsworth \& Brooks/Cole, 1990. \\
\ni [DurLM] C. Durastanti, X. Lan and D. Marinucci, Gaussian semi-parametric estimation on the unit sphere.  {\sl Bernoulli} {\bf 20} (2014), 28-77. \\
\ni [ErdMOT] A. Erd\'elyi, W. Magnus, F. Oberhettinger and F. G. Tricomi, {\sl Higher transcendental functions}, Vol, I-III.  Krieger, Melbourne, 1981 (1st ed., McGraw-Hill, 1953 (Vol. I, II), 1955 (Vol. III)). \\   
\ni [FukOT] Fukushima, M., Oshima, Y. and Takeda, M., {\sl Dirichlet forms and symmetric Markov processes}, 2nd ed., W. de Gruyter, 2011 (1st ed. 1994). \\
\ni [Gar] Garsia, A. M., Continuity properties of Gaussian processes with multidimensional time parameter. {\sl Proc. Sixth Berkeley Symp. Math. Stat. Probab. Volume II: Probability Theory} (ed. L. LeCam, J. Neyman and E. M. Scott), 369-374, U. California Press, 1972. \\
\ni [GemH] D. Geman and J. Horowitz, Occupation densities.  {\sl Ann. Prob.} {\bf 10} (1980), 1-67. \\
\ni [GinN] E. Gin\'e and R. Nickl, {\sl Mathematical foundations of infinite-dimensional statistical models}.  Cambridge University Press, 2016. \\
\ni [HauP] O. Haupt and C. Pauc, La topologie approximative de Denjoy envisag\'ee comme vraie topologie.  {\sl C. R. Acad. Sci. Paris} {\bf 234} (1952), 390-392. \\
\ni [Hor] L. H\"ormander, {\sl The analysis of linear partial differential opertors.  III.  Pseudodifferential operators}.  Grundl. math. Wiss. {\bf 274}, Springer, 1985 (corr. repr. 1994). \\
\ni [Jon] R. H. Jones, Stochastic processes on a sphere.  {\sl Ann. Math. Stat.} {\bf 34} (1963), 213-218. \\
\ni [KagLR] A. M. Kagan, Yu. V. Linnik and C. R. Rao, {\sl Characterization problems in mathematical statistics}.  Wiley, 1973. \\ 
\ni [KerOPP] Kerkyacharian, G., Ogawa, S., Petrushev, P. and Picard, D.: Regularity of Gaussian processes on Dirichlet spaces.  {\sl Constr. Approx.} {\bf 47} (2018), 277-320 (arXiv:1508.00822).\\       
\ni [Kor] J. Korevaar, {\sl Tauberian theory: A century of developments}.  Grundl. math. Wiss. {\bf 329}, Springer, 2004. \\
\ni [KraR] M. A. Krasnoselskii and Ya. B. Rutickii, {\sl Convex functions and Orlicz spaces}.  Noordhoff, Groningen, 1961 (Russian, GITTL. Moscow, 1958). \\
\ni [LanMX] Xiaohong Lan, D. Marinucci and Yimin Xiao, Strong local non-determinism and exact modulus of continuity for spherical Gaussian fields.  {\sl Stoch. Proc. Appl.} {\bf 128} (2018), 1294-1315. \\
\ni [LanX] X. Lan and Y. Xiao, Strong local non-determinism of spherical fractional Brownian motion.  {\sl Stat. Prob. Letters} {\bf 135} (2018), 44-50. \\
\ni [LangS] A. Lang and C. Schwab, Isotropic Gaussian random fields on the sphere: regularity, fast simulation and stochastic partial differential equations.  {\sl Ann. appl. Prob.} {\bf 25} (2015), 3047-3094. \\
\ni [LeoTT] N. N. Leonenko, M. S. Taqqu and G. H. Terdik, Estimation of the covariance function of Gaussian isotropic random fields on spheres, related Rosenblatt-type distributions and the cosmic variance problem.  {\sl Electronic J. Stat.} {\bf 12} (2018), 3114-3146. \\ 
\ni [Lev] P. L\'evy, {\sl Processus stochastiques et mouvement brownian}, 2nd ed., Gauthier-Villars, 1965 (1st ed. 1948). \\
\ni [LusP] H. Luschgy and G. Pag\`es, Expansions for Gaussian processes and Parseval frames.  {\sl Electronic J. Prob.} {\bf 14}, 1198-1221. \\
\ni [Ma1] Chunsheng Ma, Stationary and isotropic vector random fields on spheres.  {\sl Math. Geosciences} {\bf 44} (2012), 765-778. \\
\ni [Ma2] Chunsheng Ma, Time-varying isotropic vector random fields on spheres.  {\sl J. Theoretical Probability} {\bf 30} (2017), 1763-1785. \\
\ni [Mal1] A. A. Malyarenko, Abelian and Tauberian theorems for random fields on two-point homogeneous spaces.  {\sl Th. Prob. Math. Stat.} {\bf 69} (2005), 115-127.   \\
\ni [Mal2] A. A. Malyarenko, {\sl Invariant random fields on spaces with a group action}.  Springer, 2013.   \\
\ni [MarR] Marcus, M. B. and Rosen, J., {\sl Markov processes, Gaussian processes and local times}.  Cambridge University Press, 2006. \\
\ni [MarS1] Marcus, M. B. and Shepp, L. A.: Continuity of Gaussian processes.  {\sl Trans. Amer. Math. Soc.} {\bf 151} (1970), 377-392. \\
\ni [MarS2] Marcus, M. B. and Shepp, L. A., Sample behaviour of Gaussian processes.  {\sl Proc. Sixth Berkeley Symp. Math. Stat. Probab. Volume II: Probability Theory} (ed. L. LeCam, J. Neyman and E. M. Scott), 423-441, U. California Press, 1972. \\
\ni [MariP] D. Marinucci and G. Peccati, Random fields on the sphere.  Representation, limit theorems and cosmological applications.  {\sl London Math. Soc. Lecture Note Series} {\bf 389}, Cambridge University Press, 2011.\\
\ni [Rud] W. Rudin, {\sl Real and complex analysis}, 3rd ed.  Mc-Graw-Hill, 1987 (1st ed. 1966, 2nd ed. 1974). \\
\ni [Sak] S. Saks, {\sl Theory of the integral}, 2nd ed., Dover, 1964 (translated from Monografie Matematyczne {\bf VII}, 1937, 1st ed. Mono. Mat. {\bf II}, 1933). \\
\ni [Sch] I. J. Schoenberg, Positive definite functions on spheres.  {\sl Duke Math. J.} {\bf 9} (1942), 96-108 (reprinted in {\sl Selected Papers} (ed. C. de Boor), Birkh\"auser, 1988, Vol. 1, 172-184). \\
\ni [Sma] C. V. Smallwood, Approximate upper and lower limits.  {\sl J. Math. Anal. Appl.} {\bf 37} (1972), 223-227. \\
\ni [Ste] E. M. Stein, {\sl Singular integrals and differentiability properties of functions}.  Princeton University Press, 1970. \\
\ni [SteW] E. M. Stein and G. Weiss, {\sl Introduction to Fourier analysis on Euclidean spaces}.  Princeton Universit Press, 1971. \\
\ni [Sze] G. Szeg\H{o}, {\sl Orthogonal polynomials}.  AMS Colloquium Publications {\bf XXIII}, Amer. Math. Soc., 1934. \\                   
\ni [Yad] M. I. Yadrenko, {\sl Spectral theory of random fields}.  Optimization Software Inc., 1983. \\   
\ni [Yon] C.-H. Yong, {\sl Asymptotic behaviour of trigonometric series with modified monotone coefficients}.  Chinese University of Hong Kong, 1974. \\
\ni [Zyg] A. Zygmund, {\sl Trigonometric series}, 2nd ed., Vol. I, II.  Cambridge University Press, 1968. \\

\ni N. H. Bingham, Mathematics Department, Imperial College, London SW7 2AZ, UK; n.bingham@ic.ac.uk \\
\ni Tasmin L. Symons, Mathematics Department, Imperial College, London SW7 2AZ, UK; tls111@ic.ac.uk \\

\end{document}